\documentclass[12pt]{article}

\usepackage{amsmath}
\usepackage{amssymb}
\usepackage{amsthm}
\usepackage{hyperref} 
\usepackage[margin=1in]{geometry} 

\theoremstyle{theorem}
\newtheorem{theorem}{Theorem}
\newtheorem{lemma}{Lemma}
\newtheorem{proposition}{Proposition}
\newtheorem{corollary}{Corollary}

\theoremstyle{definition}
\newtheorem{definition}{Definition}
\newtheorem{remark}{Remark}

\newcommand{\uniform}[1]{#1\text{-uniform}}

\begin{document}

\title{The Largest Ground Set for Sunflower-Free Systems}
\author{Kyle Wood}
\date{\today} 

\maketitle

\begin{abstract}
We construct the largest $r$-flower-free system of $k$-uniform sets for $k=3, r=4$ where every element of the ground set $G$ is used at least once. The goal is to maximize $|G|$, also called the span. We define $\mathrm{sp}(k, r)$ to be the maximum $|G|$ for such a system. We show that $\mathrm{sp}(3, 4) = 39$, which yields an improved upper bound $f(3, 4) \leq 130$, improving upon the Erd\H{o}s--Rado bound of $f(3, 4) \leq 162$.
\end{abstract}

\section{Introduction}

The sunflower problem asks: what is the maximum size of a \uniform{k} set system that contains no $r$-flower? We denote this maximum by $f(k, r)$. An $r$-flower (or $r$-sunflower) is a collection of $r$ sets whose pairwise intersections are all identical; this common intersection is called the kernel.

In this paper, we introduce a new parameter $\mathrm{sp}(k, r)$, defined as the maximum size of the ground set for an $r$-flower-free system of \uniform{k} sets where every element appears in at least one set. We establish the following results:

\begin{theorem}\label{thm:main}
$\mathrm{sp}(3, 4) = 39$.
\end{theorem}

\begin{corollary}\label{cor:bound}
$f(3, 4) \leq 130$.
\end{corollary}

\section{Preliminaries}

We introduce several new concepts that will be central to our construction.

\begin{definition}
A system $\mathcal{S}$ of sets is called \emph{spiny} if each set in $\mathcal{S}$ contains at least one element that appears in no other set of $\mathcal{S}$. Such an element is called a \emph{quill}.
\end{definition}

\begin{remark}
The terminology is inspired by the hedgehog structure in hypergraphs, where each edge has a unique vertex extending outward like a spine.
\end{remark}

\begin{lemma}\label{lem:removal}
Any $r$-flower-free system remains $r$-flower-free after removing any number of sets.
\end{lemma}

\begin{proof}
Immediate from the definition of $r$-flower.
\end{proof}

\begin{lemma}\label{lem:spiny}
Any $r$-flower-free system can be reduced to a spiny system with the same span by repeatedly removing sets whose elements all appear elsewhere.
\end{lemma}

\begin{proof}
Consider a system $\mathcal{S}$ that is $r$-flower-free. While there exists a set $A \in \mathcal{S}$ such that every element of $A$ appears in at least one other set of $\mathcal{S}$, remove $A$ from $\mathcal{S}$. By Lemma~\ref{lem:removal}, the resulting system remains $r$-flower-free. 

Crucially, this process preserves the span: we only remove sets whose elements all appear elsewhere, so no element of the ground set becomes unused. The process terminates with a spiny system having the same span as the original.
\end{proof}

\begin{remark}\label{rem:spiny-importance}
Lemma~\ref{lem:spiny} is fundamental to our approach: to determine $\mathrm{sp}(k, r)$, it suffices to find the maximum span among spiny $r$-flower-free systems, since any $r$-flower-free system can be reduced to a spiny one with the same span. Thus, constructing an optimal spiny system establishes the maximum span for all systems, spiny or not.
\end{remark}

\begin{definition}
A system $\mathcal{S}$ is \emph{maximal} if no set can be added without creating an $r$-flower. A system is \emph{maximum} if it uses every element in a ground set of size $\mathrm{sp}(k, r)$.
\end{definition}

\begin{remark}
For a maximum system, each set contains exactly one quill after reduction to spiny form. If a set had multiple quills, we could add more sets using one of its quills, making it no longer a quill.
\end{remark}

\begin{lemma}\label{lem:maximal}
A maximal $r$-flower-free system must contain at least one $(r-1)$-flower with empty kernel.
\end{lemma}

\begin{proof}
If no such $(r-1)$-flower exists, we can always add a set consisting entirely of new elements, contradicting maximality.
\end{proof}

\subsection{Constructing spiny systems}

Lemma~\ref{lem:maximal} suggests a systematic approach to constructing maximal spiny systems. For concreteness, consider the case $k=3, r=4$. We must start with a 3-flower with empty kernel, say:
\[
\{1, 2, 3\}, \{4, 5, 6\}, \{7, 8, 9\}
\]

To make this system spiny while maintaining maximality, we need to designate certain elements as quills. Without loss of generality (due to symmetry), we can choose elements 7, 8, 9 to be quills, yielding:
\[
\{1, 2, 7\}, \{3, 4, 8\}, \{5, 6, 9\}
\]

This choice determines the structure of all additional sets in the system:
\begin{itemize}
\item Every additional set must contain at least one element from $\{1, 2, \ldots, 6\}$ (otherwise it forms a 4-flower with our starting 3-flower)
\item Every additional set must exclude all elements from $\{7, 8, 9\}$ (to maintain the spiny property)
\item Every additional set must contain at least one new element as its quill
\end{itemize}

This systematic construction ensures maximality while preserving the spiny structure. We can make a stronger structural change, though; once we remove the quills from such a construction, we obtain a collection of \uniform{(k-1)} sets, leading to our dimension reduction technique.

\section{The Construction Method}

\subsection{Reducing dimensions through quill removal}

The key observation is that quills in a spiny system do not contribute to intersections between sets. This allows us to establish a connection between \uniform{k} systems and \uniform{(k-1)} collections.

\begin{lemma}\label{lem:quill-removal}
Let $\mathcal{S}$ be a spiny \uniform{k} $r$-flower-free system where each set has exactly one quill. Removing all quills yields a collection $\mathcal{C}$ of \uniform{(k-1)} sets that is $r$-flower-free.
\end{lemma}

\begin{proof}
Since quills appear in exactly one set, they cannot be part of any intersection. Thus, the flower structure is preserved after their removal. Note that $\mathcal{C}$ is a collection rather than a system because removing quills may create duplicate sets.
\end{proof}

\begin{remark}\label{rem:collection-vs-system}
The distinction between collections and systems is crucial. In the standard sunflower problem, we use systems (sets are distinct) rather than collections (sets may repeat) for a fundamental reason: in a collection, we could have $r-1$ copies of each set, forming $(r-1)$-flowers whose kernel is the entire set. This would make the problem trivial. Similarly, for our construction, the optimal collection $\mathcal{C}$ should have exactly $r-1$ copies of each distinct set.
\end{remark}

\begin{lemma}\label{lem:optimal-collection}
For a maximum spiny \uniform{k} $r$-flower-free system $\mathcal{S}$, the collection $\mathcal{C}$ obtained by removing quills has each distinct \uniform{(k-1)} set appearing exactly $r-1$ times.
\end{lemma}

\begin{proof}
Each occurrence of a \uniform{(k-1)} set in $\mathcal{C}$ corresponds to a set in $\mathcal{S}$ with a unique quill. Having $r-1$ copies allows us to add $r-1$ distinct quills (one per copy), maximizing the contribution to the span. Having $r$ or more copies would create an $r$-flower in $\mathcal{C}$, while having fewer than $r-1$ copies fails to maximize the span.
\end{proof}

\subsection{The inverse construction and aura}

We can now describe the construction process in the forward direction. Given a \uniform{(k-1)} $r$-flower-free system $\mathcal{S}'$:

\begin{enumerate}
\item Create a collection $\mathcal{C}$ by taking $r-1$ copies of each set in $\mathcal{S}'$.
\item Add a unique quill to each set in $\mathcal{C}$ to obtain a spiny \uniform{k} system $\mathcal{S}$.
\end{enumerate}

The resulting span is $|G'| + |\mathcal{C}| = |G'| + |\mathcal{S}'|(r-1)$, where $G'$ is the ground set of $\mathcal{S}'$.

\begin{definition}
The \emph{aura} of a \uniform{k} system $\mathcal{S}'$ with ground set $G'$ (relative to parameter $r$) is defined as 
\[
\text{aura}(\mathcal{S}') = |G'| + |\mathcal{S}'|(r-1).
\]
\end{definition}

\begin{proposition}\label{prop:aura-characterization}
$\mathrm{sp}(k, r)$ equals the maximum aura among all \uniform{(k-1)} $r$-flower-free systems.
\end{proposition}

\begin{proof}
By Lemma~\ref{lem:spiny}, any \uniform{k} $r$-flower-free system can be reduced to a spiny one with the same span. By Lemmas~\ref{lem:quill-removal} and \ref{lem:optimal-collection}, this spiny system corresponds to a \uniform{(k-1)} system whose aura equals the span. Conversely, any \uniform{(k-1)} $r$-flower-free system can be lifted to a \uniform{k} spiny system with span equal to its aura. The maximum is therefore achieved by the \uniform{(k-1)} system with maximum aura.
\end{proof}

\begin{remark}\label{rem:aura}
Finding $\mathrm{sp}(k, r)$ requires finding the \uniform{(k-1)} $r$-flower-free system with maximum aura. This system does not necessarily have either maximum size (i.e., may not achieve $f(k-1, r)$) or maximum span (i.e., may not achieve $\mathrm{sp}(k-1, r)$).
\end{remark}

\subsection{The special case of $k=2$}

For \uniform{2} systems, there is a unique simplification: the system achieving maximum span also achieves maximum size, and therefore maximizes the aura function.

\begin{lemma}\label{lem:k2special}
For computing $\mathrm{sp}(2, r)$, the optimal \uniform{1} system $\mathcal{S}'$ simultaneously achieves both $f(1, r) = r-1$ and $\mathrm{sp}(1, r) = r-1$.
\end{lemma}

\begin{proof}
For \uniform{1} sets (singletons), any $r$-flower-free system has at most $r-1$ sets, and these sets must all be distinct singletons. Thus the unique (up to relabeling) maximum system is $\{\{1\}, \{2\}, \ldots, \{r-1\}\}$, which has both size $r-1$ and span $r-1$, yielding maximum possible aura $(r-1) + (r-1)(r-1) = r(r-1)$.
\end{proof}

\subsection{Example: Computing $\mathrm{sp}(2, 4)$}

By Lemma~\ref{lem:k2special}, we use $\mathcal{S}' = \{\{1\}, \{2\}, \{3\}\}$ with aura $3 + 3 \cdot 3 = 12$.

To construct a \uniform{2} 4-flower-free system with maximum span:

\begin{enumerate}
\item Transform $\mathcal{S}'$ to a collection $\mathcal{C}$ by taking $r-1 = 3$ copies of each set:
\[
\mathcal{C} = \{\{1\}, \{1\}, \{1\}, \{2\}, \{2\}, \{2\}, \{3\}, \{3\}, \{3\}\}
\]

\item Add a unique quill to each set in the collection:
\[
\mathcal{S} = \{\{1, 4\}, \{1, 5\}, \{1, 6\}, \{2, 7\}, \{2, 8\}, \{2, 9\}, \{3, 10\}, \{3, 11\}, \{3, 12\}\}
\]
\end{enumerate}

Therefore, $\mathrm{sp}(2, 4) = 12$.

\section{Computing $\mathrm{sp}(3, 4)$}

\subsection{Finding the optimal 2-uniform system}

Unlike the $k=2$ case, for $k=3$ we cannot assume the optimal \uniform{2} system achieves either $f(2, 4)$ or $\mathrm{sp}(2, 4)$.

From the literature, we know that $f(2, 4) = 10$. We have shown that $\mathrm{sp}(2, 4) = 12$ with a system of 9 sets. We need to find which \uniform{2} 4-flower-free system has maximum aura.

\begin{lemma}\label{lem:aura-analysis}
The \uniform{2} 4-flower-free system with maximum aura has aura 39, achieved uniquely by the system with 9 sets and span 12 and a variation of the system found by Abbott and Hanson \cite{abbott1972}.
\end{lemma}

\begin{proof}
The system achieving $\mathrm{sp}(2, 4) = 12$ has 9 sets, giving aura $12 + 9 \cdot 3 = 39$.
Abbott and Hanson \cite{abbott1972} constructed a system achieving $f(2, 4) = 10$ with span 8, giving aura $8 + 10 \cdot 3 = 38 < 39$. The only way to beat this would be increasing the span while keeping the size fixed, because with one less size, the span can be maximized as described above. 

The problem may be converted to graph theoretic terms for $k = 2$, as described by Abbott and Hanson \cite{abbott1972}. The span is equivalent to the number of vertices, and the size equivalent to the number of edges. For $r = 4$, there must be no 4-matching or vertex with degree 4 or greater. 

To tie the aura of the maximum spanning system, one would need a graph of size 10 and span 9 (10 edges, 9 non-isolated vertices). To beat the aura, one would need span 10 or greater. We use nauty \cite{mckay2014} to generate a set of all 12 vertex 10 edge graphs with maximum degree of 3. We then use a simple Python script to test whether any of these have a maximum matching less than 4. Additionally, we check how many vertices are contained in at least one edge, effectively removing isolated vertices from the graph. 

Running the script yields only a single graph with 9 vertices, and none with 10, 11 or 12 vertices (after removing isolated vertices).

Thus, only two graphs maximize the aura for $k = 2$, $r = 4$: The one which maximizes the span, and the graph found by Abbott and Hanson but with the 3-cycle component replaced by a 3-star. 
\end{proof}

\subsection{The construction}

Following our construction method from Section 3 with the optimal \uniform{2} systems having aura 39, we obtain $\mathrm{sp}(3, 4) = 39$.

\section{Application to Upper Bounds}

\begin{proof}[Proof of Corollary~\ref{cor:bound}]
Since $f(2, 4) = 10$, each element can appear in at most 10 sets of a \uniform{3} 4-flower-free system. (If an element appeared in 11 or more sets, taking all sets containing the element and removing it from each would yield a \uniform{2} system with more than 10 sets, contradicting $f(2, 4) = 10$.)

With ground set size at most $\mathrm{sp}(3, 4) = 39$ and each element appearing at most 10 times, the total number of element occurrences is at most $39 \cdot 10 = 390$. Since each \uniform{3} set contributes 3 element occurrences:
\[
f(3, 4) \leq \left\lfloor \frac{390}{3} \right\rfloor = 130.
\]
\end{proof}

This improves upon the Erd\H{o}s--Rado bound \cite{erdos1960} of $k!(r-1)^k = 6 \cdot 27 = 162$ for $k=3, r=4$.

\section{Concluding Remarks}

The uniqueness of the maximum span system for $\mathrm{sp}(3, 4)$ suggests computational approaches could further refine the upper bound. If no \uniform{3} 4-flower-free system with span 39 has more than 38 sets (verifiable computationally), then $f(3, 4) \leq \lfloor(38 \cdot 10)/3\rfloor = 126$. Because the spiny system with span 39 and the variation of Abbott and Hanson's graph are unique, any system with span 39 must contain at least one of them (up to relabeling), which is a strong restriction.

More generally, define $a(k, r)$ as the aura of a \uniform{k} $r$-flower-free system with maximum aura. Determining $a(2, r)$ would yield $\mathrm{sp}(3, r)$ and improved bounds on $f(3, r)$. Since $f(2, r)$ and $\mathrm{sp}(2, r)$ are both known, $a(2, r)$ may be quite approachable.

There are no relevant financial or non-financial competing interests to report.

\end{document}